\documentclass{amsart}[11pt]

\usepackage{amscd,amssymb,latexsym,url,verbatim,graphicx,color}
\usepackage{hyperref}
\hypersetup{colorlinks=true,linkcolor=[rgb]{0,0,1},filecolor=[rgb]{0,0,1}}
\usepackage[capitalize,nameinlink,noabbrev]{cleveref}
\usepackage{tikz,tikz-cd}

\usetikzlibrary{decorations.pathmorphing}

\usepackage{cases,amsmath}

\usepackage{pxfonts}
\usepackage[euler-digits,small]{eulervm}

\usepackage{tikz,tikz-cd}

\usepackage[dvips]{epsfig}

\title[Generalized anchored configuration spaces of graphs]
{Homology and Euler characteristic of generalized anchored configuration spaces of graphs}

\author{Dmitry N. Kozlov}

\address{Department of Mathematics, University of Bremen, 28334
  Bremen, Federal Republic of Germany.}
	
\email{dfk@math.uni-bremen.de}

\address{Okinawa Institute of Science and Technology Graduate University,
1919-1 Tancha, Onna-son, Kunigami-gun,
Okinawa, Japan.}


\newtheorem{theorem}{Theorem}[section]
\newtheorem{df}[theorem]{Definition}
\newtheorem{thm}[theorem]{Theorem} 

\newtheorem{prop}[theorem]{Proposition}
\newtheorem{lm}[theorem]{Lemma} 
\newtheorem{crl}[theorem]{Corollary}

\newtheorem{rem}[theorem]{Remark}
 
\newtheorem{fact}[theorem]{Fact} 

\newcommand{\nin}{\noindent} 
\newcommand{\pr}{\nin{\bf Proof.} }

\newcommand{\bo}{\partial}

\newcommand{\zz}{{\mathbb Z}}

\newcommand{\es}{\emptyset}

\newcommand{\ra}{\rightarrow}
\newcommand{\sm}{\setminus}
\newcommand{\supp}{\text{\rm supp}}

\newcommand{\ti}{\tilde}
\newcommand{\wti}{\widetilde}

\newcommand{\stn}[2]{\genfrac{\{}{\}}{0pt}{}{#1}{#2}}
\newcommand{\ost}[2]{#2!\stn{#1}{#2}}
\newcommand{\xor}{\,\textrm{XOR}\,}

\newcommand{\myset}[2]{\{#1\,|\,#2\}}

\graphicspath{{./figs/}}

\numberwithin{equation}{section}
\numberwithin{figure}{section}
\numberwithin{table}{section}

\newcommand{\acs}[3]{\Sigma(#1,#2,#3)}
\newcommand{\gacs}[4]{\Sigma(#1,#2,#3,#4)}

\newcommand{\cc}{\mathcal{C}}

\begin{document}

\begin{abstract}
In this paper we consider the generalized anchored configuration spaces on $n$ labeled points on a~graph.
These are the spaces of all configurations of $n$ points on a~fixed graph $G$, subject to the condition 
that at least $q$ vertices in some pre-determined set $K$ of vertices of $G$ are included in each configuration.

We give a non-alternating formula for the Euler characteristic of such spaces
for arbitrary connected graphs, which are not trees. 
Furthermore, we completely determine the homology groups of the generalized anchored 
configuration spaces of $n$ points on a circle graph.
\end{abstract}

\maketitle

\section{Introduction}

The study of the anchored configuration spaces was initiated in
\cite{Ko21} and continued in \cite{HM,Ko22}. These spaces are motivated by certain 
considerations in logistics and differ from classical configuration spaces in 
a crucial way. The formal definition is as follows.

\begin{df} 
\label{df:acs}
Let $X$ be a non-empty topological space, let $K$ be a set of distinct
points in~$X$, and let $n$ be an arbitrary positive integer.  
An {\bf anchored configuration space}, denoted $\acs XKn$,
is defined as the subspace of the direct product $X^n$, consisting of
all tuples $(x_1,\dots,x_n)$, such that $K\subseteq\{x_1,\dots,x_n\}$.
\end{df}

In this paper we continue with this line of research and generalize \cref{df:acs}, 
by relaxing the conditions on the allowed $n$-tuples $(x_1,\dots,x_n)$.

\begin{df} 
\label{df:gacs}
As above, let $X$ be a non-empty topological space, let $K$ be a set of $k$ 
distinct points in~$X$, $k\geq 0$, and let $n$ be an arbitrary positive integer. 
Furthermore, let $q$ be an arbitrary integer, such that $k\geq q\geq 0$. 
A~{\bf generalized anchored configuration space}, denoted $\gacs XKnq$,
is defined as the subspace of the direct product $X^n$, consisting of
all tuples $(x_1,\dots,x_n)$, such that $|K\cap\{x_1,\dots,x_n\}|\geq q$.
\end{df}

Clearly, the case $k=q$ in \cref{df:gacs} corresponds to \cref{df:acs}.

So far, the anchored configuration spaces have been studied in the situation when $X$ 
is a~geometric realization of a graph $G$, and $K$ is a~subset of the set of the vertices 
of~$G$. These spaces are of particular interest for 
the logistics, since they formalize networks with moving goods, with the extra condition 
that at each point in time a~certain set of nodes is securely supplied with the goods.
Accordingly, the generalized anchored configuration spaces relax this condition and only
require that at each point in time at least $q$ nodes out of some previously fixed
set $K$ are supplied.

The case of the anchored spaces, where the graph $G$ is a tree has been settled 
in \cite{Ko21}, where the homotopy type of $\acs XKn$ has been completely determined. 
Accordingly, in this paper, we consider the case when $G$ is not a~tree.

As a first result we provide a non-recursive formula for the Euler characteristic of $\gacs GKnq$, 
expressed as a sum of positive terms (rather than a sign-alternating sum). After that we turn 
to study the topology of these spaces.

We consider the case when $G$ is a circle graph, which appears to be the most natural
next step, after the case of $G$ being a~tree. This is the same as to consider the case when 
the topological space $X$ is homeomorphic to a~circle, since changing the positions 
of the points in $K$ will produce homeomorphic anchored configuration spaces. 
Since all we need to record is the cardinality of $K$, we let $\Omega(k,n)$ denote
$\acs GKn$, where $G$ is a cycle graph with $k$ vertices, and $K$ is the set of vertices of~$G$.
Passing on to the generalized anchored configuration spaces, we let $\Omega(k,n,q)$ denote 
$\gacs GKnq$ in this case.

The spaces $\Omega(2,n)$ were the focus of investigations in \cite{Ko22} and \cite{HM}.
More specifically, the homology of these spaces was calculated in \cite{Ko22} using 
discrete Morse theory. This work was continued in \cite{HM}, where the cup product structure
was completely described, and connection to the topological complexity was established.

In this paper we study the spaces $\Omega(k,n)$ of an arbitrary $k$, and more generally
$\Omega(k,n,q)$, for an arbitrary $q\leq k$, and calculate their homology groups in all 
dimensions. Rather than using discrete Morse theory,  our method is to consider 
classical long sequences for the corresponding combinatorially given chain complexes. 
For the standard concepts of Algebraic Topology we refer to \cite{Fu,GH,Hat,Mu}. 
Our study lies within the field of Applied Topology, see \cite{Ca,EH,Ko08,Ko20} 
for more information.



\section{The Euler characteristic of the generalized anchored configuration spaces}

Let $G$ be a~connected graph, which is not a~tree, and let $V$ and $E$ denote its sets 
of vertices and edges, respectively. Let $K$ be an arbitrary subset of $V$, and let $q\leq |K|$.
In this section we give a~theorem which provides a complete non-recursive and
non-alternating formula for the Euler characteristic of the spaces $\gacs GKnq$. 
Before we proceed with its formulation and its proof, let recall the following concepts.

First, for arbitrary positive integers $a\geq b$ the {\it Stirling numbers of 
the second kind}, denoted $\stn{a}{b}$, count the number of ways to partition 
a set of $a$ labelled objects into $b$ nonempty unlabelled subsets.
Clearly, then $b!\stn{a}{b}$ is the number of ways to partition 
a set of $a$ labelled objects into $b$ nonempty labelled subsets.

Second, if we have a set $U$, a subset $S\subseteq U$ and an element $x\in U$,
we let $S\xor x$ denote the subset of $U$ obtained from $S$ by the {\it exclusive or}
operation with respect to $x$. Formally, we set
\[
S\xor x:=
\begin{cases}
S\sm x,& \text{ if }x\in S; \\
S\cup x,& \text{ if }x\not\in S.
\end{cases}
\]


We can now formulate the main result of this section.

\begin{thm} \label{thm:ecmain}
Let $G$ be an arbitrary connected graph,  whose set of vertices is $V$, 
and whose set of edges is $E$. 
Let $K$ be an arbitrary non-empty subset of $V$, and let $q$ be a positive integer, such that $q\leq |K|$. 
Finally, let $n$ be a natural number, such that $n\geq q$. 

Assume the graph $G$ is not a tree. Then, the Euler characteristic of the cell complex $\gacs GKnq$
is given by the formula\footnote{Note, that in the formula \eqref{eq:genec} we use the convention $0^0=1$, 
while of course $0^t=0$, for $t>0$.} \label{p:0}
\begin{equation}
\label{eq:genec}
\frac{(-1)^{n-q}}{q!}\chi(\gacs GKnq)=
\sum_{\lambda=\varepsilon}^{\varepsilon+k-q}\sum_{t=0}^{n-q}\binom{\lambda-\varepsilon+q-1}{q-1}
\binom{n}{t}\stn{n-t}{q}\lambda^t,
\end{equation}
where $k:=|K|$ and $\varepsilon:=|E|-|V|$.
\end{thm}
\pr 
Set $\alpha:=|V|$, $\beta:=|E|$, and let us write $V=\{v_1,\dots,v_\alpha\}$, $E=\{e_1,\dots,e_\beta\}$.
Since $G$ is a connected graph which is not a~tree, we have $\beta\geq\alpha$, or,
using our notation $\varepsilon=\beta-\alpha\geq 0$. Without the loss of generality
we can assume that $K=\{v_1,\dots,v_k\}$.

We can think about the Euler characteristic $\chi(\gacs GKnq)$ as a sum of
$\pm 1$, more precisely $(-1)^{\dim\sigma}$, ranging over the set of all cells $\sigma$
in  $\gacs GKnq$. By definition of $\gacs GKnq$, each such cell is indexed by a~function 
$\varphi:V\cup E\ra 2^{[n]}$, which satisfies two conditions: 
\begin{enumerate}
\item[(1)] the number of vertices $v\in K$, for which $\varphi(v)\neq\emptyset$ is at least $q$;
\item[(2)] the set of images $\myset{\varphi(x)}{x\in V\cup E}$ is a~partition of 
$[n]=\{1,\dots,n\}$ into disjoint sets.
\end{enumerate}

Let us now introduce some further notation. 
Consider the following collection of sets:
\begin{align*}
A_i:=&\varphi(v_i)\cup\varphi(e_i),\textrm{ for }1\leq i\leq \alpha, \\
U:=&\bigcup_{i=\alpha+1}^\beta \varphi(e_i).
\end{align*}
Set $P_\sigma:=(A_1,\dots,A_\alpha,U)$. Clearly, 
the tuple $P_\sigma$ is an ordered set partition of $[n]$,
in which we allow empty sets.
We shall now group all the cells $\sigma\in\gacs GKnq$ according to their
tuple $P_\sigma$, and calculate the contribution to the Euler characteristic
separately in each group.

Consider first an arbitrary tuple $P_\sigma$, such that $\cup_{i=k+1}^\alpha A_i\neq\emptyset$. 
Let $M$ be the set of all cells with this tuple. Let $l$ be the minimal element of $\cup_{i=k+1}^\alpha A_i$,
and let $t$ denote the index $k+1\leq t\leq \alpha$, for which $l\in A_t=\varphi(v_t)\cup\varphi(e_t)$. 
We can then define an involution $\mu:M\ra M$, by moving the element $l$ from 
$\varphi(v_t)$ to $\varphi(e_t)$, and vice versa.
Formally, we set 
\[\mu(\varphi)(u):=
\begin{cases}
\varphi(u)\xor l,&\text{ if } u=v_t, \text{ or } u=e_t;\\
\varphi(u), & \text{ otherwise}. 
\end{cases}\]
Since $k+1\leq t\leq\alpha$, there are no conditions on $\varphi(v_t)$, so the involution $\mu$ 
is well-defined. It produces a perfect matching on the set $M$. The difference of dimensions of 
any two matched cells is $1$, so their contributions to the Euler characteristic
of $\gacs GKnq$ have opposite signs. It follows that the contribution of each matched
pair is $0$, and hence also the total contribution of all the cells in $M$ is $0$.

This means, that when computing the Euler characteristic of $\gacs GKnq$ we can limit
ourselves to considering the tuples $P_\sigma$, for which $A_{k+1}=\dots=A_{\alpha}=\emptyset$,
which we do for the rest of the argument.


Assume now $\sigma$ is one of the remaining cells. Set 
\[\{i_1,\dots,i_m\}:=\myset{1\leq i\leq k}{\varphi(v_i)\neq\es},\text{ where }i_1<\dots<i_m,\] 
and set $r(\sigma):=i_q$. This is well-defined since by condition (1) above, we know that $m\geq q$. 

Let us now fix the following data $\Omega$: 
\begin{itemize}
\item the index set $\{i_1,\dots,i_q\}$, where $i_1<\dots<i_q$, 
\item the non-empty sets $A_{i_1}$, $\dots$, $A_{i_q}$.
\end{itemize}
Let $M$ denote the set of all cells with this data $\Omega$ (and with $A_{k+1}=\dots=A_\alpha=\es$).
Note, that for each cell $\sigma\in M$, we have $r(\sigma)=i_q$.  
Let us calculate the total contribution of the cells in $M$ to the Euler characteristic of $\gacs GKnq$. 

Let $\wti M$ denote the subset of $M$ consisting of all cells $\sigma$ for which
the union $\bigcup_{j=r(\sigma)+1}^s A_j$ is not empty. For $\sigma\in\wti M$, set
$\rho(\sigma):=\min\bigcup_{j={r(\sigma)+1}}^s A_j$. Let $s\geq t>r(\sigma)$ be the index, 
for which $\rho(\sigma)\in A_t$. In a complete analogy to the above, we define a matching
$\mu:\wti M\ra\wti M$ by moving the element $\rho(\sigma)$ from $\varphi(v_t)$ to $\varphi(e_t)$ 
and vice versa. This is a perfect matching on $\wti M$, since $t\geq i_q+1$, so there is no 
restriction $\varphi(v_t)$ being non-empty.  The difference of the dimensions of the matched 
cells is equal to $1$. This implies that the total contribution to the Euler characteristic by the cells 
from $\wti M$ is $0$. We can therefore from now on concentrate on the cells from $M\sm\wti M$.

For $1\leq j\leq q$, we set $l_j$ to be the minimum of $A_{i_j}$. Since
$A_{i_j}\neq\es$, the element $l_j$ is well-defined.

We now partition the set $M\sm\wti M$ into the sets $M_1,M_2,\dots,M_{q+1}$ as follows.
For each cell $\sigma\in M\sm\wti M$ we define $h(\sigma)$ to be the index $z$, uniquely
determined by the following condition:
\[\varphi(v_{i_z})\neq l_z,\text{ and }\varphi(v_{i_j})=l_j,\text{ for all }j<z.\]
Here, if $\varphi(v_{i_j})=l_j$, for all $1\leq j\leq k$, we set $h(\sigma)=q+1$.
Clearly $1\leq h(\sigma)\leq q+1$, and we define the above partition of $M\sm\wti M$ 
by saying that  $\sigma\in M_i$ if and only if $h(\sigma)=i$.

Next, fix an index $1\leq d\leq q$, and calculate the contribution of the cells in $M_d$.
Same way as earlier in the proof, we can define an involution $\mu:M_d\ra M_d$. This time
it is shifting $l_d$ from $\varphi(v_{i_d})$ to $\varphi(e_{i_d})$ and back. Formally,
\[\mu(\varphi)(u):=
\begin{cases}
\varphi(u)\xor l_d,&\text{ if } u=v_{i_d}, \text{ or } u=e_{i_d};\\
\varphi(u), & \text{ otherwise}. 
\end{cases}\]
Since $\varphi(v_{i_d})\neq l_d$, the involution $\mu$ is well-defined. As before,
it matches cells with dimension difference $1$, so the contribution of these two cells,
and hence also the contribution of the total set $M_d$ to the Euler characteristic of 
$\gacs GKnq$ is~$0$.

The only interesting contribution occurs in $M_{q+1}$. Note, that all cells in $M_{q+1}$
have dimension $n-q$. Indeed, if $\sigma\in M_{q+1}$, we have $\varphi(v_{i_j})=l_j$, 
for all $1\leq j\leq q$, and $\varphi(v_t)=\emptyset$, for $t\notin\{i_1,\dots,i_q\}$.
It follows that  $\sum_{v\in V}|\varphi(v)|=q$, hence $\sum_{e\in E}|\varphi(e)|=n-q$.

This means that each $\sigma\in M_{q+1}$ gives the contribution $(-1)^{n-q}$, and we
need to compute the cardinality $|M_{q+1}|$. Set $W:=[n]\sm\bigcup_{j=1}^q A_{i_j}$.
The cells $\sigma\in M_{q+1}$ are obtained by arbitrarily distributing the elements of $W$
among the sets $\varphi(e_j)$, for
\begin{itemize}
\item either $j\in\{1,\dots,r(\sigma)\}\sm\{i_1,\dots,i_q\}$, 
\item or $\alpha+1\leq j\leq \beta$.
\end{itemize}
In total, there are $\beta-\alpha+r(\sigma)-q=\varepsilon+r(\sigma)-q$ such sets, so
we have $|M_{q+1}|=(\varepsilon+r(\sigma)-q)^{|W|}$. 

At this point, let us specifically consider what happens when $\varepsilon+r(\sigma)-q=0$,
which of course is equivalent to saying that $\varepsilon=0$ and $r(\sigma)=q$. In this case, 
there are no sets to distribute the elements of $W$ to. Therefore, the number of ways to 
distribute the elements of $W$, and hence also the cardinality of $M_{q+1}$, is equal to~$0$,
unless, of course, the set $W$ itself is empty, in which case the cardinality of $M_{q+1}$ is
equal to~$1$. Note, how this is compatible with our convention for $0^t$, cf. the footnote 
on page~\pageref{p:0}.

Summing over all choices of $\Omega$, we have
\begin{equation} \label{eq:o2}
\chi(\gacs GKnq)=\sum_{\Omega}(-1)^{n-q}(\varepsilon+i_q-q)^{|W|}.
\end{equation}

To further evaluate \eqref{eq:o2} we can choose the data $\Omega$ in the following order. 
First, pick $r$, such that $q\leq r\leq k$. Set $i_q:=r$ and choose the remaining elements 
$i_1,\dots,i_{q-1}$, such that $i_1<\dots<i_{q-1}<i_q$, in $\binom{r-1}{q-1}$ ways. After that, 
choose the cardinality $t:=|W|$, we have $0\leq t\leq n-q$. Proceed by choosing $W$ itself, 
there are $\binom{n}{t}$ possibilities. Finally, distribute the elements of $[n]\sm W$ into 
the sets $A_{i_1},\dots,A_{i_q}$, so that they are non-empty. The number of ways to do that is
$q!\stn{n-t}{q}$.

Summarizing, we obtain
\[\chi(\gacs GKnq)=
(-1)^{n-q}\sum_{r=q}^{k}\binom{r-1}{q-1}\sum_{t=0}^{n-q}
\binom{n}{t}q!\stn{n-t}{q}(\varepsilon+r-q)^t.\]
Now, set $\lambda:=\varepsilon+r-q$. Then $r=q,\dots,k$ translates to 
$\lambda=\varepsilon,\dots,\varepsilon+k-q$, and $r-1=\lambda-\varepsilon+q-1$, so we obtain 
\eqref{eq:genec}.
\qed

We can now specialize \cref{thm:ecmain} to the case of the regular anchored configuration spaces.

\begin{crl}
The Euler characteristic of $\acs GKn$ is given by the formula
\begin{multline}
\label{eq:ecnt}
\frac{(-1)^{n-k}}{k!}\chi(\acs GKn)=\sum_{t=0}^{n-k}\binom{n}{t}\stn{n-t}{k}\varepsilon^t=\\
\stn nk+\binom{n}{1}\stn{n-1}{k}\varepsilon+\binom{n}{2}\stn{n-2}{k}\varepsilon^2+
\cdots+\binom{n}{n-k}\stn kk\varepsilon^{n-k}.
\end{multline}
\end{crl}
\pr
Substitute $q:=k$ into \eqref{eq:genec}. We then have $\lambda=\varepsilon$, so the first summation
is trivial, and the $\binom{\lambda-\varepsilon+q-1}{q-1}=\binom{q-1}{q-1}=1$. This yields \eqref{eq:ecnt}.
\qed

\section{The chain complexes for the generalized anchored configuration spaces on circle graphs}

Let us fix positive integers $k$ and $n$, such that $n\geq k\geq 2$. 
Let $C_k$ be a cycle graph with $k$ vertices and $k$ edges.
Let $E$ denote its set of edges, and let $V$ denote its set of vertices.
We can choose the index set to be $\zz_k$, and write $E=\{e_1,\dots,e_k\}$ and 
$V=\{v_1,\dots,v_k\}$, in such a~way that the adjacency map $\bo:E\ra 2^V$ is 
given by $\bo(e_i)=\{v_i,v_{i+1}\}$.\footnote{Of course, $k+1=1$ in $\zz_k$.} 

\begin{df} \label{df:nt}
Given $n$, a {\bf vertex-edge $n$-tuple} is an $n$-tuple 
$\sigma=(\sigma_1,\dots,\sigma_n)$, such that $\sigma_i\in V\cup E$, for all $i$.

For a vertex-edge $n$-tuple $\sigma=(\sigma_1,\dots,\sigma_n)$ we define two subsets 
of $\zz_k$, which we call {\bf vertex } and {\bf edge support sets}, and which we denote 
$\supp_V(\sigma)$ and $\supp_E(\sigma)$, as follows:
\[\supp_V(\sigma):=\myset{i\in\zz_k}{v_i\in\{\sigma_1,\dots,\sigma_n\}},\]
and
\[\supp_E(\sigma):=\myset{j\in\zz_k}{e_j\in\{\sigma_1,\dots,\sigma_n\}}.\]

Finally, the {\bf dimension} of $\sigma$ is defined to be $\dim\sigma:=|\myset{i}{\sigma_i\in E}|$.
So, in particular, we have $0\leq\dim\sigma\leq n$.
\end{df}

Clearly, for any vertex-set $n$-tuple $\sigma=(\sigma_1,\dots,\sigma_n)$, 
the set $\{\sigma_1,\dots,\sigma_n\}$ is a disjoint union of the sets 
$\myset{v_i}{i\in\supp_V(\sigma)}$ and $\myset{e_j}{j\in\supp_E(\sigma)}$.

The direct product $\underbrace{C_k\times\dots\times C_k}_n$ has a~natural structure 
of the cubical complex, whose geometric realization is an $n$-torus. 
Its cells are indexed by the vertex-edge $n$-tuples, whose dimensions, as described 
in \cref{df:nt}, coincide with the geometric dimension of the corresponding cells. 
Therefore, the chain complex whose chain groups are generated by the vertex-edge 
$n$-tuples, with appropriately defined boundary operators, will calculate the homology 
of an $n$-torus.

We shall now consider the chain complexes whose chain groups are generated  
by the vertex-edge $n$-tuples, satisfying additional conditions on the vertex support
set $\supp_V(\sigma)$.

\begin{df}
Assume we are given an arbitrary subset $P\subseteq\zz_k$, and a~nonnegative integer $q$, 
such that $q\leq |P|$.
We define a chain complex $\cc^{P,q}=(C^{P,q}_*,\bo_*)$, where $C^{P,q}_*$ are 
free abelian groups, as follows.

\begin{enumerate}
\item[(1)] For each $d$, the free abelian group $C^{P,q}_d$ is generated by
the vertex-edge $n$-tuples $\sigma=(\sigma_1,\dots,\sigma_n)$, 
with $\dim\sigma=d$, satisfying the following two conditions:
\begin{itemize}
\item $\supp_V(\sigma)\subseteq P$;
\item $|\supp_V(\sigma)|\geq q$.
\end{itemize}

\item[(2)]
The boundary operator takes the vertex-edge $n$-tuple $\sigma$, and replaces, 
with an appropriate sign,
any of the edges $\sigma_i\in E$ by any of its boundary vertices, subject to 
the condition that the index of that vertex lies in~$P$. Formally we have
\begin{equation} \label{eq:bo1}
\bo\sigma=\sum_{i\in\supp_E(\sigma)}\sum_{\ti\sigma_i\in\bo\sigma_i\cap V_P}
(-1)^{\rho(\sigma,i)}
(\sigma_1,\dots,\sigma_{i-1},\tilde\sigma_i,\sigma_{i+1},\dots,\sigma_n),
\end{equation}
where $V_P:=\myset{v_j}{j\in P}$, and 
$\rho(\sigma,i):=|\supp_E(\sigma)\cap\{1,\dots,i-1\}|$.
\end{enumerate}
\end{df}

Note the special case when $|P|=q$, when the chain groups $C_d^{P,q}$ are 
generated by all $\sigma$, satisfying $\dim\sigma=d$ and $\supp_V(\sigma)=P$.

For convenience, we introduce additional notation for the complement set 
$H:=\zz_k\sm P$, and $h:=|H|=n-|P|$.

\begin{rem} \label{rem:star}
Obviously, $C_i^{P,q}=0$, for $i<0$. Furthermore,  if a~vertex-edge $n$-tuple 
 $\sigma=(\sigma_1,\dots,\sigma_n)$ satisfies $|\supp_V(\sigma)|\geq q$, then 
$\dim\sigma\leq n-q$, so $C_i^{P,q}=0$ also for all $i>n-q$.
\end{rem}

In what follows, we shall compute the homology groups of the chain complexes $\cc^{P,q}$.
When $P$ is a proper subset of $\zz_k$, the complexes $\cc^{P,q}$ do not correspond
to topological spaces, and play here an auxilliary role. Accordingly, the case which interests 
us most is when $P=\zz_k$, since it gives us the homology of the generalized anchored 
configuration spaces $\Omega(k,n,q)$. We stress this observation for a~later reference.
\begin{fact} \label{fact}
The chain complex $\cc^{\zz_k,q}$ is isomorphic to the cubical chain complex
of the generalized anchored configuration space $\Omega(k,n,q)$.
In particular, the chain complex $\cc^{\zz_k,k}$ is isomorphic to the cubical chain complex
of the anchored configuration space $\Omega(k,n)$.
\end{fact}
Our calculation will proceed by induction, and we shall compute the homology groups 
for all values of $P$ and~$q$.


\section{Calculation of the homology groups of $\cc^{P,q}$}

\subsection{The case $q=0$} $\,$

\nin
Let us start with the case $q=0$. When $q=0$ the condition $|\supp_V(\sigma)|\geq q$ 
is void, which radically simplifies the situation. The homology is then given by 
the following proposition.

\begin{prop} \label{prop:q0}
 $\,$

\begin{enumerate}
\item[(1)] The chain complex $\cc^{\zz_k,0}$ calculates the homology of an $n$-torus. 
In fact, it is a chain complex of the cubical complex obtained as an $n$-fold direct product 
of the $k$-cycle.
\item[(2)] When $P$ is a proper subset of $\zz_k$, we have
$H_n(\cc^{P,0})\approx\zz^{h^n}$, and all other homology groups are trivial.
\end{enumerate}
\end{prop}
\pr Statement $(1)$ is trivial and simply formalizes our earlier observation, so we proceed 
to proving the statement $(2)$. 

Let $\wti C_k$ be the graph which is in a sense dual to $C_k$. It is also a cycle graph 
with $k$ vertices and $k$ edges, but with a different indexing. Let $\wti E$ 
denote its set of edges, and let $\wti V$ denote its set of vertices. Both again are 
indexed by $\zz_k$, $\wti E=\{\ti e_1,\dots,\ti e_k\}$ and $\wti V=\{\ti v_1,\dots,\ti v_k\}$, 
but now in such a~way, that the boundary map $\bo:\wti E\ra 2^{\wti V}$ is 
given by $\bo(\ti e_i)=\{\ti v_{i-1},\ti v_i\}$. So, compared to $C_k$, the relative indexing 
is shifted by~$1$.

Let $G$ denote the subgraph of the cycle graph $\wti C_k$, obtained by deleting all edges 
indexed by $H$. Consider the cubical complex $G^n=\underbrace{G\times\dots\times G}_{n}$, 
and consider the cochain complex of $G^n$, let us call it $\tilde\cc^*$. 
It is easy to see that $\cc^{P,0}$ is isomorphic to this cochain complex, with 
the isomorphism $\varphi$ given by $\varphi(v_i):=\ti e_i$, and
$\varphi(e_i):=\ti v_i$, for all $i\in\zz_k$. In particular, we have
$H_i(\cc^{P,0})\approx H^{n-i}(\ti\cc^*)$, for all~$i$.

On the other hand, we assumed that $h\geq 1$, so topologically, the graph $G$ consists 
of $h$ disjoint intervals. In particular, the direct product $G^n$ is homotopy equivalent to 
the discrete space with $h^n$ points. Therefore, we have
\[H^i(\ti\cc^*)=\begin{cases}
\zz^{h^n},&\text{ if } i=0;\\
0, &\text{ otherwise,} 
\end{cases}
\]
and it follows that 
\[H_i(\cc^{P,0})=\begin{cases}
\zz^{h^n},&\text{ if } i=n;\\
0, &\text{ otherwise.} 
\end{cases}
\]
\qed

\subsection{Structure of the relative chain complexes} $\,$

\nin
Assume now $q\geq 1$, and consider the chain complex $\cc^{P,q-1}$. 
The condition as to which vertex-edge $n$-tuples are allowed to be taken as generators
of the chain groups is weaker for $\cc^{P,q-1}$, than it is for $\cc^{P,q}$, so 
the latter is its chain subcomplex. The following lemma states that their quotient 
can be decomposed as a direct sum of chain complexes of the same type.

\begin{lm} \label{lm:quot}
For any $P\subseteq\zz_k$, and any $q\geq 1$, we have the following chain complex isomorphism:
\begin{equation} \label{eq:quot}
\cc^{P,q-1}/\cc^{P,q}\approx\bigoplus_S\cc^{S,q-1},
\end{equation} 
where the sum is taken over all subsets of $P$ of cardinality $q-1$.
\end{lm}
\pr For each $d$, the relative chain group 
$C_d(\cc^{P,q-1}/\cc^{P,q})=C_d^{P,q-1}/C_d^{P,q}$ is generated by the cosets of
$C_d^{P,q}$, whose representatives are the vertex-edge $n$-tuples 
$\sigma=(\sigma_1,\dots,\sigma_n)$, satisfying $|\supp_V(\sigma)|=q-1$.

Call such a coset $\bar\sigma$. The relative boundary operator in 
$\cc^{P,q-1}/\cc^{P,q}$ is then given by the following formula, cf.\ \eqref{eq:bo1},

\begin{equation} \label{eq:bo2}
\bo\bar\sigma=\sum_{i\in\supp_E(\sigma)}\sum_{\ti\sigma_i\in\bo\sigma_i\cap \bar V}
(-1)^{\rho(\sigma,i)}
(\sigma_1,\dots,\sigma_{i-1},\tilde\sigma_i,\sigma_{i+1},\dots,\sigma_n),
\end{equation}
where $\bar V=\myset{v_j}{j\in\supp_V(\sigma)}$, and $\rho(\sigma,i)$ is the same as 
in \eqref{eq:bo1}.

In other words, when taking the boundary, we are allowed to replace an edge
with any of its boundary vertices, subject to the condition, that this does
not change the vertex support set.

Since the boundary operator preserves the vertex support set, the chain complex
$\cc^{P,q-1}/\cc^{P,q}$ decomposes as a direct sum, with direct summands indexed
by all possible choices of $\supp_V(\sigma)$, which is the same as to say all
possible choices of $(q-1)$-subsets of~$P$. This proves~\eqref{eq:quot}.
\qed 


\subsection{The case $P\neq\zz_k$} $\,$

\nin
When $P$ is a proper subset of $\zz_k$, it turns out that all the homology
of the chain complex $\cc^{P,q}$ is concentrated in its top dimension.

\begin{thm} \label{thm:h1}
Assume $P$ is a proper subset of $\zz_k$. 
Then, the homology of $\cc^{P,q}$ is concentrated in dimension $n-q$,
in other words, $H_i(\cc^{P,q})=0$, for $i\neq n-q$.
\end{thm}
\pr The proof proceeds by induction on $q$. For the base case $q=0$, 
this has been proved in \cref{prop:q0}(2).

Assume now $q\geq 1$. Since the chain complex $\cc^{P,q}$ is a~subcomplex 
of $\cc^{P,q-1}$, we have the following long exact sequence:
\begin{equation}
\label{eq:1}
\dots\ra H_*(\cc^{P,q})\ra H_*(\cc^{P,q-1})\ra H_*(\cc^{P,q-1}/\cc^{P,q})
\stackrel{\partial}{\ra} H_{*-1}(\cc^{P,q})\ra\dots
\end{equation}

Note, that by induction assumption, the homology of the complex $\cc^{P,q-1}$ is concentrated 
in dimension $n-(q-1)=n-q+1$. Furthermore, due to dimensional reasons, see \cref{rem:star}, 
the homology of $\cc^{P,q}$ must be $0$ in dimension $n-q+1$ and above.

By \cref{lm:quot} we have $\cc^{P,q-1}/\cc^{P,q}\approx\oplus_S\cc^{S,q-1}$, 
where the sum is taken over all subsets of $P$ of cardinality $q-1$.
Since each $S$ is a proper subset of $\zz_k$, by induction assumption, the homology of 
$\cc^{S,q-1}$ is also concentrated in dimension $n-q+1$.
It follows that the only nontrivial part of the long exact sequence \eqref{eq:1} is
\[0\ra H_{n-q+1}(\cc^{P,q-1})\ra H_{n-q+1}(\cc^{P,q-1}/\cc^{P,q})\ra H_{n-q}(\cc^{P,q})\ra 0,\]
so it follows that $H_i(\cc^{P,q})=0$, for $i\neq n-q$.
\qed


\subsection{The case $P=\zz_k$} $\,$

\nin We are now ready to deal with the main case.

\begin{thm} \label{thm:main}

The homology of the chain complex $\cc^{\zz_k,q}$ is given by 
the following formula:
\begin{equation} \label{eq:mainf}
H_i(\cc^{\zz_k,q})\cong
\begin{cases}
\zz^{\binom{n}{i}},& \text{ if } 0\leq i\leq n-q-1, \\
0,& \text{ if } i<0\text{ or }i>n-q.
\end{cases}
\end{equation}
\end{thm}

\pr
Once again, we proceed by induction on $q$.
When $q=0$, we simply have the homology of the $n$-torus, see \cref{prop:q0}(1).

Assume now $q\geq 1$. Consider again the long exact sequence:
\begin{equation}
\label{eq:2}
\dots\ra H_*(\cc^{\zz_k,q})\ra H_*(\cc^{\zz_k,q-1})\ra H_*(\cc^{\zz_k,q-1}/\cc^{\zz_k,q})
\stackrel{\partial}{\ra} H_{*-1}(\cc^{\zz_k,q})\ra\dots
\end{equation}

\cref{lm:quot} together with \cref{thm:h1} imply that $H_i(\cc^{\zz_k,q-1}/\cc^{\zz_k,q})=0$,
for all $i\neq n-q+1$. Furthermore, for dimensional reasons, we have $C_i^{\zz_k,q}=0$,
whenever $i<0$, or $i>n-q$, see \cref{rem:star}, so we know that we must have 
$H_i(\cc^{\zz_k,q})=0$, unless $0\leq i\leq n-q$.

It follows that the long exact sequence \eqref{eq:2} falls into several short 
pieces $H_i(\cc^{\zz_k,q})\approx H_i(\cc^{\zz_k,q-1})$, for $0\leq i\leq n-q-1$, 
and one longer piece
\begin{multline}
0\ra H_{n-q+1}(\cc^{\zz_k,q-1})\ra H_{n-q+1}(\cc^{\zz_k,q-1}/\cc^{\zz_k,q})\ra \\
H_{n-q}(\cc^{\zz_k,q}) \ra H_{n-q}(\cc^{\zz_k,q-1})\ra 0.
\end{multline} 
This implies the statement of the theorem.
\qed

Note, that for dimensional reasons, see \cref{rem:star}, the top-dimensional 
homology group $H_{n-q}(\cc^{\zz_k,q})$ must be free.
The Betti number $\beta_{n-q}(\cc^{\zz_k,q})$ can then be
computed using the Euler-Poincar\'e formula.


When $G$ is a cycle, we have $|V|=|E|$. This means that $\varepsilon=0$, so 
in \eqref{eq:ecnt} all the terms except for the first one vanish, and we have 
the following formula.

\begin{crl}
\label{crl:euler}
We have $\chi(\Omega(n,k))=(-1)^{n-k}k!\stn{n}{k}$.
\end{crl}

Together with \cref{thm:main} this finishes the calculation of the Betti numbers of $\Omega(n,k)$.

For the generalized anchored configuration spaces, we need to substitute $\varepsilon=0$ in
\eqref{eq:genec}, and obtain the following corollary.
\begin{crl}
\label{crl:euler2}
Assume $1\leq q\leq k$, we have 
\[\chi(\Omega(n,k,q))=(-1)^{n-q}q!\sum_{\lambda=0}^{k-q}
\sum_{t=0}^{n-q}\binom{\lambda+q-1}{q-1}\binom nt\stn{n-t}{q}\lambda^t.\]
\end{crl}

Again, together with \cref{thm:main} this gives us the Betti numbers of $\Omega(n,k,q)$.

\newpage

\nin
{\bf Declarations.}

\vskip5pt 

\nin
{\bf Ethical Approval:} not applicable.

\vskip5pt 

\nin
{\bf Competing interests:} not applicable.

\vskip5pt 

\nin
{\bf Authors' contributions:} not applicable.

\vskip5pt 

\nin
{\bf Funding:} DFG Project Nr. 509120879

\vskip5pt 

\nin
{\bf Availability of data and materials:} not applicable.


\end{document}